\documentclass[11pt]{article}
\usepackage{amssymb}
\usepackage{mathrsfs}
\usepackage{amsfonts}
\usepackage{tikz}

\usepackage{amsmath,amsthm,amsfonts,amssymb,amscd}

\allowdisplaybreaks[4]
\newtheorem {Problem} {Problem}[section]
\newtheorem {Theorem} [Problem]{Theorem}
\newtheorem {Lemma}[Problem]{Lemma}
\newtheorem{Conjecture}[Problem]{Conjecture}

\newtheorem {Corollary}[Problem]{Corollary}
\newenvironment {Proof}{\noindent {\bf Proof.}}{\hfill\ensuremath{\square}}
\newcommand*{\QEDB}{\hfill\ensuremath{\square}}

\begin{document}

\title{ The spectral radius of  minor free graphs \thanks{This work is supported by  Hainan Provincial Natural Science Foundation of China (Nos. 120RC453, 120MS002) and the National Natural Science Foundation of China (Nos. 11971311, 12026230), the
Montenegrin-Chinese Science and Technology Cooperation Project (No.3-12).
}}

\author{ Ming-Zhu Chen, A-Ming Liu,\\
School of Science, Hainan University, Haikou 570228, P. R. China, \\
\and  Xiao-Dong Zhang\footnote{Corresponding author. E-mail: xiaodong@sjtu.edu.cn (X.-D. Zhang)}
\\School of Mathematical Sciences, MOE-LSC, SHL-MAC\\
Shanghai Jiao Tong University,
Shanghai 200240, P. R. China\\
Email: chenmingzhuabc@163.com, aming8809@163.com, xiaodong@sjtu.edu.cn.\\
\\
}
\date{}
\maketitle

\begin{abstract}
In this paper,  we present a sharp upper bound for the spectral radius of an $n$-vertex  graph without $F$-minor for sufficient large $n$, where $F$ is obtained from the complete graph $K_r$ by deleting disjointed paths. Furthermore, the graphs which achieved the sharp bound are characterized. This result may be regarded to be an extended revision of the number of edges in an $n$-vertex graph without $F$-minor.


\end{abstract}

\section{Introduction}

 The study of $H$-minor free graphs has occupied a center place in  graph theory.
It is a well-known cornerstone result that  a graph is planar if and only if it is $K_5$-minor free and   $K_{3,3}$-minor free \cite{Wagner}. Moreover, a  graph is outerplanar if and only if it is $K_4$-minor free and $K_{2,3}$-minor free \cite{BM}.

In extremal graph theory,  center problems ask to maximize or minimize a graph invariant over a fixed family of graphs.
 One of the well-known problem is to determine the maximum number of edges among $H$-minor free graphs and characterize all such graphs attaining this maximum number of edges.  There is an key result on $H$-minor free graphs
 \begin{Theorem}\cite{Mader}
 \label{maderth}
Let $G$ be any  $K_r$-minor free graph on $n$ vertex. Then the number $e(G)$ of edges in $G$ at most
  $(r-2)n-\binom{r-1}{2}$  for $2\le r\le 7$ and $n\ge r$.
  \end{Theorem}
 For $r\le 5,$ this was first proved by Dirac [5]. Some years later but independently of Mader,
Gy\"{o}ri \cite{gyori} proved Theorem~\ref{maderth} 1.1 for $r=6$. Furthermore, Mader \cite{Mader} pointed out Theorem~\ref{maderth} does not hold for $r\ge 8$.   Later, Seymour and Thomas (see \cite{songT2006}) conjectured the following:
\begin{Conjecture}
 For every $r\ge 1$,  there exists a constant $N$ such that every $(r-2)$
connected  and $K_r$-minor free  graph on $n\ge N$  vertices  has  at most  $(r-2)n-\binom{r-1}{2}$ edges.
\end{Conjecture}

Let $K_r^-$  and $K_r^{*}$ denote the graph obtained from $K_r$  by deleting one edge and any two edges, respectively.
As a weaker variant of Seymour-Thomas's conjecture,  the maximum number of edges in graphs without  $K_r^-$-minor or  $K_r^*$ has
been studied.  Dirac \cite{dirac} proved  the following for $r=5$ and $r=6$.  Jakoblen \cite{jakobsen1972, jakobsen1983} proved the cases for $r=7$   and $r=8$ with $K_8^*$. Song \cite{song2005} proved the case for $r=8$ and $K_8^-$.
\begin{Theorem}\cite{dirac, jakobsen1972, jakobsen1983, song2005}For $r\in \{5, \ldots,8\}$, if $G$ is an $n\ge r-1$-vertex graph with
 at least  $\frac{1}{2}[(2r-5)n-(r-3)(r-1)]$  $($or $\frac{1}{2}[(2r-6)n-(r-4)(r-1)]$ $)$ edges, then $G$ contains $K_r^-$minor ($K_r^*$-minor), except some known graphs.
\end{Theorem}
Recently,  Rolek \cite{rolek2019} proved the above result still holds for $r=9$ and $K_9^*$-minor. Up to now, the maximum edges of graphs without  $K_r^-$minor (or $K_r^*$-minor) for $r\ge9$  (or $r\ge 10$ respectively) still open.

In spectral extremal graph theory,  Hong \cite{Hong} presented a sharp upper bound for the  spectral radius of all $K_5$-minor free graphs and  characterized all extremal graphs. Fang \cite{Fang} obtained upper bound for the spectral radius of graphs without $K_{3,3}$-minor. Later,
Nikiforov \cite{Nikiforov} gave a sharp bound for the spectral radius of graphs without $K_{2,t}$-minor and characterized all extremal graphs which achieved the bound.
Recently,  Tait \cite{Tait2019} proved the two key results on  $K_r$-minor and $K_{s,t}$-minor, which extended  and strengthened the above results.
 \begin{Theorem}\cite{Tait2019}
 Let $r\ge 3.$  There exist an integer $N$ such that the  $n(\ge N)$-vertex graph with no $K_r$ minor of maximum spectral radius is the join of $K_{r-2}$ an independent set of size $n-r+2$.
 \end{Theorem}

  This paper is motivated by the maximum number of edges of graphs without $K_r^-$-minor in extremal graph theory and the maximum spectral radius of graphs without $K_r$-minor. We will discuss a related question. If $G$ is $n$-vertex graph without $K_r^-$-minor, how large can the spectral radius of its adjacency matrix be ?   For more information on relationships between the spectral radius and minor free, the readers may be referred to
\cite{chen2019+, DNP1,Hong, Nikiforov10, Tait2019}.

 Let $G$ be an undirected simple graph with vertex set
$V(G)=\{v_1,\dots,v_n\}$ and edge set $E(G)$.
The \emph{adjacency matrix}
$A(G)=(a_{ij})$ of $G$  is the $n\times n$ matrix, where
$a_{ij}=1$ if $v_i$ is adjacent to $v_j$, and $0$ otherwise. The  \emph{spectral radius} of $G$ is the largest eigenvalue of $A(G)$.  For $v\in V(G)$,  the \emph{neighborhood} $N_G(v)$ of $v$  is $\{u: uv\in E(G)\}$ and the \emph{degree} $d_G(v)$ of $v$  is $|N_G(v)|$.
We write $N(v)$ and $d_v$ for $N_G(v)$ and $d_G(v)$ respectively if there is no ambiguity.
A graph $H$ is a \emph{minor} of a graph $G$ if $H$ can be obtained from $G$ by a sequence of vertex
and edge deletions and edge contractions. A graph $G$ is \emph{$H$-minor free} if   it  does not contain $H$ as a minor.
  For two vertex disjoint graphs $G$ and $H$,  we denote by  $G\cup H$ and  $G\nabla H$  the \emph{union} of $G$ and $H$,
and the \emph{join} of $G$ and $H$ which is obtained by joining every vertex of $G$ to every vertex of $H$, respectively.
Denote by $k G$  the $k$ disjoint union of $G$. For an edge subset $M$ of a graph $G$,
$G-M$ denotes the graph obtained from $G$ by deleting the edges in $M$,
while for an edge subset $M'$ of the complement of $G$,
$G+M'$ denotes the graph obtained from $G$ by adding the edges in $M'$.  Denote $S_{n,t}:=K_{t}\nabla \overline{K}_{n-t}$, i.e.,  $S_{n,t}$ is the join of $K_t$ and an independent set of size $n-t$; and  denote $F_{n,t}:=K_{t}\nabla (p K_2\cup K_q)$, where $t\geq3$, $n-t=2p+q$ and $0\leq q<1$.
 In addition, denote  $K_r^{=}$  the graph  obtained from the complete graph $K_r$ by deleting any two incident edges.  Note that $K_r^{=}$ is different from $K_r^*$ which is the graph from the complete graph $K_r$ by deleting any two  edges (maybe independent edges or incident edges).
  For graph notation and terminology undefined here,  readers are referred to \cite{BM}.

The main result of this paper can be stated as follows.

%
%
\begin{Theorem}\label{thm1}
Let $F$ be a graph obtained from the complete graph $K_r$ by deleting the disjoint paths $P_{k_i}$ on $k_i$ vertex for $1\le i\le l$, i.e., $F=K_r-E(\cup_{i=1}^l P_{k_i})$, where $V(P_i)\bigcap V(P_j)=\emptyset$ for $1\le i\neq j\le l$.
 Suppose that $G$ is an $F$-minor free graph of  sufficiently large order  $n$.

(i) If $ k_1=\cdots=k_l=2$ , then
$$\lambda(G)\leq \lambda(F_{n,r-3})$$
 with equality if and only if $G=F_{n,r-3}$.

(ii) If there is at least one $k_i\geq3$, then
$$\lambda(G)\leq
  \lambda(S_{n,r-3})$$ with equality if and only if $G=S_{n,r-3}$.
\end{Theorem}

{\bf Remark} It is easy to see (for example see \cite{CLZ2019} or \cite{Nikiforov10}) that
$$\lambda(S_{n,r-3})=\frac{r-4+\sqrt{4(r-3)n-(3r^2-16r+20)}}{2}.$$
In addition, if $n-r+3$ is even, then $$\lambda(F_{n,r-3})=\frac{r-4+\sqrt{4(r-3)n-(3r^2-14r+11)}}{2};$$
and if $n-r+3$ is odd, then $\lambda(F_{n,r-3})$ is the largest root of

$$x^3-(r-3)x^3-[(r-3)n-(r^2-5r+5)]x+r-3=0.$$

The rest of this paper is organized as follows. In Section 2, some preliminary are presented. In Sections  3 and 4, we give a proof of Theorem~\ref{thm1} (i) and (ii) respectively. In Section 5, we conclude by observing that Theorem~\ref{thm1} can be furthermore extended.

\section{Preliminary}


Mader \cite{Mader} proved an elegant result on the number of edges in $H$-minor graphs.
\begin{Lemma}\label{e-minor1}\cite{Mader}
Let $G$ be an $n$-vertex graph. For every graph $H$, if $G$ is $H$-minor free,  then   there exists a constant $C$  such that $$e(G)\leq Cn.$$
\end{Lemma}



\begin{Lemma}\label{e-minor3}
Let $G$ be an $n$-vertex $K_r^-$-minor  (or $K_r^{=}$-minor) free bipartite graph with vertex partition $A$ and $B$. If $|A|=k$ and $|B|=n-k$, then there exists a constant $C$ depending only on $r$ such that
$$e(G)\leq Ck+(r-3)n.$$
\end{Lemma}


\begin{Proof}
Since a $K_r^{=}$-minor free bipartite graph is  also $K_r^{-}$-minor free, it suffices to prove that the result holds for the $K_r^{-}$-minor free bipartite graphs. So we suppose that $G$ is  $K_r^{-}$-minor free. This yields that $G$ is  $K_{r-2,r}$-minor free as
 any $K_{r-2,r}$-minor contains  a $K_r^{-}$-minor. Hence the assertion follows from   \cite[Lemma~2.1]{Tait2019} or \cite[Theorem 2.2]{Thomason2007}. 
  \end{Proof}

\begin{Lemma}\label{spec1-minor3}
Let  $G$ be an $n$-vertex  graph having  the maximum spectral radius $\lambda(G)$  among all  $K_r^-$minor (or $K_r^{=}$-minor) free connected graphs with $n\ge r\ge  4$. Then

(i) $\lambda(G)\geq \sqrt{(r-3)(n-r+3)}$.

(ii) If  $\mathbf x=(x_u)_{u\in V(G)}$ is a positive eigenvector with the maximum entry $1$ which corresponds to   $\lambda(G)$,  then  $x_u\geq \frac{1}{\lambda (G)}$ for all $u\in V(G)$.
\end{Lemma}

\begin{Proof}
Since a $K_r^{=}$-minor free graph is  also $K_r^{-}$-minor free, it suffices to prove that the assertion   holds for  $K_r^{-}$-minor free graphs.  So we suppose that $G$ is  $K_r^{-}$-minor free.

 (i) Since $K_{r-3,n-r+3}$ is  $K_r^-$-minor free and $G$ has the maximum  spectral radius $\lambda(G)$  among all  $K_r^-$minor (or $K_r^{=}$-minor) free connected graphs, we have
 $$\lambda(G)\geq \lambda(K_{r-3,n-r+3})\geq \sqrt{(r-3)(n-r+3)}.$$

(ii)  Let $w\in V(G)$ such that $x_w=1$.
If $u=w$, then $x_u=1\geq \frac{1}{\lambda (G)}$. So next we suppose that $u\neq w$. We consider the following two cases.

{\bf Case~1.} $u$ is adjacent to $w$. By the eigenequation of $G$ on $u$, $$\lambda(G)x_u=\sum_{uv\in E(G)}x_v\geq x_w=1,$$
which implies that $x_u\geq \frac{1}{\lambda (G)}.$

{\bf Case~2.} $u$ is not adjacent to $w$. Let $G'$ be the graph obtained from $G$ by deleting all edges incident with $u$ and adding an edge $uw$.
Note that $uw$ is a pendant edge.  If $G$ is  $K_r^-$-minor free, then  $G'$ is also $K_r^-$-minor free.
Hence
\begin{eqnarray*}
  0 &\geq& \lambda(G')-\lambda(G)\geq \frac{\mathbf x^{\mathrm{T}}A(G')\mathbf x}{\mathbf x^{\mathrm{T}}\mathbf x}- \frac{ {\bf x^T} A(G)\bf x}{\bf {x^T}\bf x}\\
   &=& \frac{2x_u}{\mathbf x^{\mathrm{T}}\mathbf x} \Big(x_w-\sum_{uv\in E(G)}x_v\Big)\\
 &=&  \frac{2}{\mathbf x^{\mathrm{T}}\mathbf x} \Big(1-\sum_{uv\in E(G)}x_v \Big),
\end{eqnarray*}
which implies that $$\sum_{uv\in E(G)}x_v\geq1.$$
By  the eigenequation of $G$ on $u$, $$\lambda(G) x_u=\sum_{uv\in E(G)}x_v\geq1,$$ which yields that $$x_u\geq \frac{1}{\lambda (G)}.$$
\end{Proof}


\section{Proof of Theorem~\ref{thm1}~(i) }
In order to prove Theorem~\ref{thm1}~(i), we first give some structural properties of $K_r^-$-minor free graphs.
\begin{Lemma}\label{L-minor=}
Let $G$ be an $n-$vertex  $K_r^{-}$-minor free graph  with $ r\geq4$. If $G$ contains a complete bipartite subgraph $K_{r-3,(1-\delta)n}=[A, B]$ with 
 $|A|=r-3$ and $|B|=(1-\delta)n>r$, then the induced subgraph by vertex set $B$ in $G$ is $P_3$-free and there are at least  $(1-3\delta)n$ vertices in $B$ which has no neighbors in $V(G)\backslash (A\cup B)$.
\end{Lemma}

\begin{Proof} We first prove that the induced subgraph $G[B]$ is $P_3$-free. In fact, suppose that there exists a path  $P_3=uvw$ in $G[B]$.
Since $[A,B]$ is the complete bipartite  and $|B|\ge |A|+3$, there exists a matching $M$ in $[A,B]$ such that $M$ saturates each vertex in $A$ and does not saturate $\{u,v,w\}$.  Vertex set $\{u,v,w\}$ and contacting all edges in $M$ in $ G$ forms a $K_r^-$, which is a contradiction. So   $G[B]$ is $P_3$-free.


Furthermore, we have the following Claim.

{\bf Claim.} If  $C$ is a  connected component of $G-(A\cup B)$. then there are at most two vertices in $B$ which have  a neighbor in $C$ respectively.

Suppose that there are three vertices $u,v,w\in B$  which have  a neighbor in $C$ respectively.  Then there exist two paths $P$ from $u$ and $v$ and
$Q$ from $v$ to $w$ such that all interior vertices of both $P$ and $Q$  are in $C$. Moreover, there exists a matching $M$ in $[A,B]$ such that $M$ saturates each vertex in $A$ and does not saturate $\{u,v,w\}$.  Then vertex set $\{u,v,w\}$ and contacting all edges in $M$ in $ G$ forms a $K_r^-$, which is a contradiction. So Claim holds.

Now  let $R=V(G)\backslash (A\cup B)$ and
 $D=\{v\in B: N_G(v)\cap R=\emptyset \}.$
 By the definition of $R$, $$|R|=n-|A|-|B|\leq n-(r-3)-(1-\delta)n<\delta n,$$ which implies that
$R$ has at most $ \delta n$ components. By Claim, $B\backslash D$ has at most $ 2\delta n$ vertices.
Hence $$|D|=|B|- |B\backslash D|\geq (1- \delta) n-2\delta n\geq(1- 3\delta) n.$$
This completes the proof.
\end{Proof}

\begin{Theorem}\label{edge-minor-1}
Let  $G$ be an $n$-vertex $K_r^{-}$-minor free graph with $ r\geq4$. If $G$ contains a complete bipartite subgraph $K_{r-3,(1-\delta)n}= [A,B]$  with $|A|=r-3$, $|B|=(1-\delta)n$ and $(1-3\delta)n>r$. Let $G^*$  be the graph  obtained from $G$ by adding edges to $A$ to make it a clique. Then
 $G^*$ is also  $K_r^-$-minor free. \end{Theorem}

\begin{Proof}
%
%
 Denote by $R=V(G)\backslash (A\cup B)$ and $D=\{v\in B: N_G(v)\cap R=\emptyset \}$. Suppose that  $G^*$ contains a  $K_r^-$-minor. Then there exist $r$  disjoint vertex sets $X_1, X_2,\ldots, X_r\subseteq V(G^*)=V(G))$ with the following properties:

  (a) $G^*[X_i]$ is connected for $i=1, \dots, r$.

  (b) There has at least one edge between each pair distinct sets  $X_i$ and $X_j$ in $G^*$ except  one pair sets.

 (c) There exists a set $X_j$ such that $X_j\bigcap D \not= \emptyset$, where  $1\le j\le r$.

 In fact, if $X_1, \ldots, X_r\subseteq (B\setminus D)\bigcup R$, then $X_1, \ldots X_r$  in $G$ form $K_r^-$, a contradiction. Hence there exists a set $X_j\bigcap A\not=\emptyset$, then choose a vertex $u\in D$ and let $X_j^{\prime}=X_j\bigcup \{u\}$. Then
 $X_1, \ldots, X_j^{\prime}, \ldots, X_r$ satisfying (a), (b) and (c).
  Furthermore, let
  $$f(X_1,\ldots,X_r)=|\{i:|X_i\bigcap D|\geq1 ~\text{for}~ 1\leq i\leq r\}|.$$
  Hence  we can choose $r$ disjoint sets  $X_1, \ldots, X_r$  satisfying  (a) and (b) such that  $f(X_1,\ldots,X_r)$ is as large as possible and
  $$|X_1\bigcap D|\ge |X_2\bigcap D|\ge \cdots\ge |X_r\bigcap D|.$$

  {\bf Claim~1.} If  $|X_m\bigcap D|\ge 2$ and $|X_{m+1}\bigcap D|\le 1$, then $X_i\bigcap A \not= \emptyset$ for $i=1, \ldots, m$.

   Suppose there exists a $X_j\bigcap A=\emptyset$ for $1\le j\le m$.
      Then $G[X_j]=G^*[X_j]$ is connected.  If there exist two adjacent $u_1$ and $u_2$ vertices in $G[X_j\bigcap D]$, then    by $|A|=r-3$, there there must exist two sets $X_s,X_t$ such that $X_i\bigcap A= \emptyset$ for $i=s,t$ and $s,t\neq j$. Furthermore, by (b), there is at lest one edge either  between $X_j$ and $X_s$ or  between $X_j$ and $X_t$. Note that there is no edge between $D$ and $R$. So $G[B]$ contains a path of order $3$, which contradicts to Lemma~\ref{L-minor=}. Hence there exist two non-adjacent vertices $u_1$ and $u_2 $ in $X_j\bigcap D$.  Then there exists a path $P$ from $u_1$ to $u_2$ in $G[X_j]$ by $G[X_j]$ being connected. Note that $X_j\subseteq B\bigcup R$ and there is no edge between $D$ and $R$ in $G$. Thus the first edge $u_1v_1$ and last edge $v_2u_2$ in $P$ are independent edges in $B$.  Moreover, by Lemma~\ref{L-minor=}, $|D|\ge (1-3\delta)n>r$. We can choose $(r-3)$  sets $\{a_i,b_i\}$ for $a_i\in A$ and $b_i\in D\setminus \{u_1,u_2\}$ for $i=1, \ldots, r-3$,  and $\{u_1\}$, $\{u_2\}$, a path $V(P)\setminus\{u_1,u_2\}$ which form a $K_r^-$ in $G$. This is a contradiction.  So Claim 1 holds.

      Furthermore, choose a vertex $w_i\in X_i\bigcap D$ for $i=1, \ldots, m$.  Let  $Y_i=(X_i\setminus (X_i\bigcap D))\bigcup\{w_i\}$, for $i=1, \ldots, m$,  and $Y_i=X_i$ for $i=m+1, \ldots, r$.  Then by Claim~1, it is easy to see  that the following Claim 2 holds.

      {\bf Claim~2.}  $Y_i\bigcap A\not= \emptyset$ for $i=1, \ldots, m$; $f(Y_1,\ldots,Y_r)=f(X_1, \ldots, X_r)$
  and $|Y_i\bigcap D|\le 1$.

      {\bf Claim~3.}  $G^*[Y_i]$ is connected for $i=1, \dots, r$.

      Since $G^*[Y_i]=G^*[X_i]$ is connected for $m+1\le i\le r$. We only consider $G^*[Y_i]$ for $1\le i\le m$. By Claim~2, $Y_i\bigcap A\not=\emptyset$. For any  two vertices $u,v$ in $Y_i$, there exists a path $P$  from $u$ to $v$  in $G^*[X_i]$ since $G[X_i]$ is connected. If $P$ contains a vertex $b\in X_i\bigcap D\setminus\{w_i\}$, then there exist two vertices $x,y$ in $X_i$ such that they are adjacent to $b$ and $x,y\notin R$. Since $G[B]$ is $P_3$-free, either $x$ or $y$ is in $A$. So $x$ is adjacent to $y$. So there is a path in $X_i$  from $u$ to $v$ containing no $b$. Hence $G[Y_i]$ is connected and Claim 3 holds.

    {\bf Claim~4.}  There has at least one edge between each pair distinct sets  $Y_i$ and $Y_j$ in $G^*$ except  one pair sets.

 If $1\le i\neq j\le m$, then by Claim~2, $Y_i\bigcap A\neq \emptyset$ and  $Y_j\bigcap A\neq \emptyset$, which implies that there has at least one edge between $Y_i$ and $Y_j$ since $G^*[A]$ is a complete graph.  If $1\le i\le m$ and $1\le j\le f(Y_1, \ldots, Y_r)$, then $Y_i\bigcap A\neq \emptyset$  and $Y_j\bigcap\neq \emptyset.$  So there has at least one edge between $Y_i$ and $Y_j$ Since $G[A,D]$ is a complete graph.
 If $m+1\le i\neq j\le r$, then   There has at least one edge between each pair distinct sets  $Y_i$ and $Y_j$ in $G^*$ if and only if  there has at least one edge between each pair distinct sets  $X_i$ and $X_j$ in $G^*$. Hence we only consider that  $1\le i\neq j\le m$ and  $f(Y_1,\ldots, Y_r)+1 \le j\le r$.  It is easy to see that if either $Y_j\bigcap A\neq \emptyset$ or $ Y_j \bigcap B\neq \emptyset $ then there has at least one edge between $Y_i$ and $Y_j$ since a vertex of $Y_i$ is adjacent to each vertex in $A\bigcup B$ in $G^*$.  Therefore we may assume that $Y_j\subseteq R$ which implied that there  has no edge between $D$ and $Y_j$.
 Hence there has at least one edge between each pair distinct sets  $Y_i$ and $Y_j$ in $G^*$ if and only if  there has at least one edge between each pair distinct sets  $X_i$ and $X_j$ in $G^*$. So Claim~4 holds.

 {\bf Claim~5.}  If $Y_i\cap A\neq \emptyset$,  then $Y_i\cap D\neq \emptyset$, for $1\leq i\leq r$.

  By Claims~ 3 and 4, disjoint sets $Y_1, \ldots, Y_r$ satisfy (a), (b) and (c). Now suppose that there exists $1\le j\le r$ such that $Y_j\cap A\neq \emptyset$ and $Y_j\cap D= \emptyset$. Then choose a vertex $w\in D$ and let $Z_j=Y_j\bigcup\{w\}$ and $Z_i=Y_i$ for $1\le i\neq j\le r$. It is easy to see that the $r$ disjoint sets $Z_1, \ldots, Z_r$ satisfy (a), (b) and (c). Moreover,
  $$f(Z_1, \ldots, Z_r)=f(Y_1, \ldots, Y_r)+1=f(X_1, \ldots, X_r)+1$$
  which contradicts to the choice of $X_1, \ldots, X_r$.  So Claim~5 holds.

   {\bf Claim~6.} $G[Y_i]$ is connected for $i=1, \ldots, r$.

  Since   $G^*[Y_i]$ is connected,  there exists a path $P$ from $u$ to $v$ in $G^*[Y_i]$ for any two vertices $u,v \in Y_i$. If the path $P$ contains an edge $a_1a_2$ with $a_1, a_2\in A$. By Claim~5,  $Y_i\bigcap D\neq \emptyset$ and there exists a vertex $w\in Y_i \bigcap D$. So the edge of the path $P$ may be replaced by edges $a_1w$ and $a_2w$, which yields a walk $Q$ from $u$ to $v$ which contains no edge in $A$. So there exists a path from $u$ to $v$ in $G[Y_i]$ and  $G[Y_i]$ is connected.

  {\bf Claim~7.}  There has at least one edge between each pair distinct sets  $Y_i$ and $Y_j$ in $G$ except  one pair sets.

By Claim~4, suppose that there has  at least one edge between $X_i$ and $X_j$ in $G^*$ and there has no edges between $X_i$ and $X_j$ in $G$. Then there must exist two vertices $y_i,y_j $ such that $y_i\in Y_i$ and $y_j\in Y_j.$ By Claim~5, there exists a vertex $w\in Y_i\bigcap D$. Hence there has one edge $y_jw$ between $Y_i$ and $Y_j$ in $G$.  So Claim~7 holds.

By Claims 6 and 7, $Y_1, \ldots, Y_r$ form a $K_r^-$minor of $G$ which is a contradiction. We finish our proof.
  \end{Proof}

\begin{Theorem}\label{edge-minor-3}
Let  $G$ be a $K_r^{-}$-minor free graph of order $n$ with $ r\geq4$.  If there exist two disjoint subsets $A,B\subseteq V(G)$ such that $G[A]$  is a complete graph and $G[A,B]$ is a complete bipartite graph, where $|A|=r-3$, $|B|\geq1$,  then
 $G^*=G-\{uv:v\in N_B(u)\}+\{uv:v\in A\setminus N_A(u)\}$  is $K_r^-$-minor free for any  $u\in V(G)\backslash  (A\cup B)$, where $N_A(u)=N_G(u)\bigcap A$ and $N_B(u)=N_G(u)\bigcap B$
 \end{Theorem}

 \begin{Proof}
Suppose that $G^*$ contains a  $K_r^-$-minor. Then there exist $r$ disjoint  vertex subsets $X_1,\ldots, X_r\subseteq V(G^*)$ with the following properties:

(a) For each $X_i$, $G^*[X_i]$ is connected.

  (b) There has at least one edge between each pair distinct sets  $X_i$ and $X_j$ in $G^*$ except  one pair sets.

It is easy to see that $u\in X_s$ for some $1\le s\le r$. Otherwise, $u\notin  \bigcup _{i=1}^r X_i$  and $X_1,\ldots, X_r$ in $G$ form a $K_r^-$ which is a contradiction.
Furthermore, by  (b),  there has at least $(r-2)$ edges in $[X_s, \bigcup_{i=1}^r\setminus X_s]$. Note that $d_G^*(u)=|A|=r-3\ge 1$. So $X_s$ has at least two vertices.  there exists a vertex $a\in A\bigcap X_s$  by $G^*[X_s]$ being connected.

{\bf Claim~1.} $G[X_i]$ is connected for $i=1, \ldots, r$.

We only prove that $G[X_s]$ is connected. For any two vertices $x,y$ of $ X_s$ in $G$. there exists a path $P$ from $x$ to $y$ in $G^*$ since $G^*[X_s]$ is connected. If $P$ in $G^*$ contains an edge $uw$ with $w\in N_G(u)\bigcap B$, then there is a walk from $x$ to $y$ in $G$ obtained from $P$ replaced $uw$ by a path $uaw$. Hence $G[X_s]$ is connected.

{\bf Claim~2.} There has at least one edge between each pair distinct sets  $X_i$ and $X_j$ in $G$ except  one pair sets.

We only prove that there has at least one edge between $X_s$ and $X_j$ in $G$ if and only if  there has at least one edge between $X_s$ and $X_j$ in $G^*$.  Without loss of generality assume that there has one edge $uw$ in $X_s$ and $X_j$, where $w\in X_j\bigcap B$. Then there has on edge $aw$ where $a\in A\bigcap X_s$ and $w\in X_j$. So Claim 2 holds.

By Claims~1 and 2, $X_1, \ldots, X_r$ in $G$ form a $K_r^-$ which is a contradiction. Hence we finish our proof.
\end{Proof}

{\bf Remark} Theorems~\ref{L-minor=} and \ref{edge-minor-3} are in their own interest. Now we are ready to prove the main results.


{\bf Proof of Theorem~\ref{thm1}~(i).}

 Let $G$ be an $F$-minor free graph of  sufficiently large order  $n$ with the maximum spectral radius.
Then  $G$ must be connected. Otherwise adding an edge to  components of $G$ to make it connected  will result in a new $F$-minor free graph with larger spectral radius, a contradiction. Now we consider the following two cases.

{\bf Case~1.} $l=1$, i.e., $F=K_r^-$. Let  and $\mathbf x=(x_u)_{u\in V(G)}$ be a positive eigenvector with  the maximum entry $x_w=1$  corresponding   to    $\lambda (G)$.  For $0< \epsilon< 1$, denote $$L=\{v\in V(G): x_v> \epsilon\}, S=\{v\in V(G): x_v\leq \epsilon\}, $$
where $\epsilon$  will be chosen later. Clearly,
$V(G)=L\cup S$.
 By Lemma~\ref{e-minor1}, there is a constant $C>r$ such that
\begin{equation}\label{1-1}
\begin{aligned}
2e(S)\leq 2e(G)\leq Cn.
   \end{aligned}
 \end{equation}
In addition, by Lemma~\ref{spec1-minor3},
\begin{equation}\label{1}
\begin{aligned}
\lambda(G)\geq \sqrt{(r-3)(n-r+3)}.
   \end{aligned}
 \end{equation}

{\bf Claim~1.} 
$e(L,S)\leq (r-3+\epsilon)n$ and $2e(L)\leq \epsilon n$.

It is easy to see that
\begin{eqnarray*}
  \lambda(G)\varepsilon |L| &< & \sum_{v\in L}\lambda(G) x_v =   \sum_{v\in L} \sum_{vz\in E(G)} x_z
   \leq   \sum_{v\in L}d_v\leq 2e(G),
\end{eqnarray*}
which implies that
\begin{equation}\label{upper L} |L|\le \frac{2e(G)}{\varepsilon\lambda(G)}\le \frac{C\sqrt{n}}{\epsilon\sqrt{ r-3.5}}
\end{equation}
by (\ref{1-1}) and (\ref{1}).
Furthermore, by Lemma~\ref{e-minor3}, there is a constant $C'$ only depending on $r$ such that
\begin{equation}\label{2-1}
\begin{aligned}
e(L,S)\leq C'|L|+(r-3)n\leq \frac{C'C\sqrt{n}}{\epsilon\sqrt{ r-3.5}}+(r-3)n\leq (r-3+\epsilon )n,
  \end{aligned}
 \end{equation}
  as long as  $\epsilon\geq \sqrt{C'C}(r-3.5)^{-\frac{1}{4}}n^{-\frac{1}{4}}$.
In addition, by (\ref{upper L}) and  Lemma~\ref{e-minor1}, we have
\begin{equation}\label{2}
\begin{aligned}
2e(L)\leq C|L|\leq \frac{C^2\sqrt{n}}{\epsilon\sqrt{ r-3.5}}\leq\epsilon n,
  \end{aligned}
 \end{equation}
 as long as $\epsilon\geq C(r-3.5)^{-\frac{1}{4}}n^{-\frac{1}{4}}$.
So Claim~1 hold.


{\bf Claim~2.}
If $u\in L$  with $x_u=1-\delta$, then  there exists a constant $C_1$  independent of  $\delta$ and $\epsilon$ such that
$$ d_u\geq (1- C_1(\delta+\epsilon))n.$$

Clearly
\begin{eqnarray*}
  \lambda \sum_{v\in V(G)}x_v &=& 
  \sum_{v\in V(G)} \sum_{vz\in E(G)}x_z
  = \sum_{v\in V(G)} d_vx_v
  \\ &\leq &\sum_{v\in L} d_v+\epsilon \sum_{v\in S} d_v\\
   &=&  2e(L)+2\epsilon e(S)+(1+\epsilon)e(L,S),
\end{eqnarray*}
which implies
\begin{equation}\label{I1}
 \sum_{v\in V(G)}x_v\leq \frac{ 2e(L)+2\epsilon e(S)+(1+\epsilon)e(L,S)}{\lambda}.
\end{equation}
Denote by $B_u=\{v\in V(G):uv\notin E(G)\}$. By Lemma~\ref{spec1-minor3}~(ii) and (\ref{I1}),
\begin{equation*}\label{I2}
  \begin{aligned}
 \frac{1}{  \lambda }|B_u|&\leq \sum_{v\in B_u}x_v
\leq \sum_{v\in V(G)}x_v-  \sum_{uv\in E(G)}x_v =\sum_{v\in V(G)}x_v-\lambda x_u\\
   &\leq  \frac{ 2e(L)+2\epsilon e(S)+(1+\epsilon)e(L,S)}{\lambda}-\lambda x_u.
     \end{aligned}
\end{equation*}
Furthermore, using  (\ref{1-1}), (\ref{2-1})  and  (\ref{2}), we have
\begin{eqnarray*}
  |B_u| &\leq&  2e(L)+2\epsilon e(S)+(1+\epsilon)e(L,S)-\lambda^2x_u \\
   &\leq& \epsilon n+\epsilon \cdot Cn+(1+\epsilon)(r-3+\epsilon)n-(r-3)(n-r+3)(1-\delta)  \\
   &=& \Big(\epsilon(1+C)+(1+\epsilon)(r-3+\epsilon)-(r-3)(1-\delta)\Big)n+(r-3)^2(1-\delta)\\
   &\leq&(C+r)(\delta+\epsilon) n+(\delta+\epsilon)n\\
   &=& (C+r+1)(\delta+\epsilon)n,
\end{eqnarray*}
where the last second inequality holds as long as  $\epsilon\geq \frac{(r-3)^2}{n}$.
Hence $$d_u=n-1-|B_u|\geq n-1-(C+r+1)(\delta+\epsilon)n\geq (1-(C+r+2)(\delta+\epsilon))n.$$
Let $C_1=C+r+2$ which is independent of $\delta$ and $\varepsilon$.  So  Claim~2 holds.


{\bf Claim~3.}
There exist $r-3$ vertices $v_1, \ldots, v_{r-3}\in L$ which satisfy
$x_{v_i}\ge 1-C_3\varepsilon$ and $d_{v_i}\ge (1-C_3\varepsilon)n$ for $1\le i\le r-3$, where $C_3$ is  a constant independent of $\varepsilon$ and $n$.

 We prove Claim 3 by the induction method. First  choose $v_1=w\in L$. Then $x_{v_1}=1$. Furthermore, by Claim 2, there exists a constant $c_1>1$  independent of $\varepsilon$ and $n$ such that
 $d_{v_1}\ge (1-c_1\varepsilon)n$.

Now assume that we have chosen  $v_1, \ldots, v_k\in L$ which  satisfy
$x_{v_i}\ge 1-c_2\varepsilon$ and $d_{v_i}\ge (1-c_2\varepsilon)n$ for $1\le i\le k$, where $c_2$  is  a constant independent of $\varepsilon$ and $n$.
Denote  $\eta=c_2\varepsilon$ and $U=\{v_1, \ldots, v_k\}.$


By the eigenvector-eigenvalue equation for $A(G)^2$,
\begin{eqnarray*}
  (r-3)(n-r+3) &\leq&\lambda^2(G)x_w=\sum\limits_{vw\in E(G)}\sum\limits_{vz\in E({G})}x_z\leq\sum\limits_{vz\in E(G)}(x_v+x_z)\\
   &=& \sum\limits_{vz\in E(S)}(x_v+x_z) +\sum\limits_{vz\in E(L,S)}(x_v+x_z)+\sum\limits_{vz\in E(L)}(x_v+x_z)\\
    & \leq& 2\epsilon e(S)+2e(L)+\epsilon e(L,S)+\sum_{\substack{ uv\in E( L\backslash U,S)\\u \in L\backslash U }}x_u+\sum_{\substack{uv\in E(L\cap U,S)\\u\in L\cap U}}x_u.\\
    &\le & 2\varepsilon cn+\varepsilon n+ \varepsilon (r-3+\varepsilon)n+kn+\sum_{\substack{ uv\in E( L\backslash U,S)\\u \in L\backslash U }}x_u,
\end{eqnarray*}
which implies that
\begin{equation}\label{claim3-2} \sum_{\substack{ uv\in E( L\backslash U,S)\\u \in L\backslash U }}x_u\ge (r-3-k-\varepsilon(C+r))n,\end{equation}
    as long as  $\epsilon \geq \frac{(r-3)^2}{n}$.
  On the other hand, since
$$e(L\cap U, S)+e(L\cap U, L\backslash U)+2e(L\cap U)=\sum_{v\in L\cap U} d(v)\geq k(1-\eta)n,$$
we have
\begin{eqnarray*}
  e(L\cap U, S)  &\geq& k(1-\eta)n-e(L\cap U,L\backslash U)-2e(L\cap U) \\
   &\geq& k(1-\eta)n-k(|L|-k)-k(k-1) \\
   &\geq&  k(1-\eta)n-k(\epsilon n-t)-k(k-1)\\
      &=& k(1-\eta-\epsilon )n+k.
\end{eqnarray*}
Hence  
 \begin{equation}\label{claim3-3} e(L\backslash U,S)= e(L,S)-e(L\cap U, S)\leq  (r-3+\epsilon)n-k(1-\eta-\epsilon )n-k <(r-3+\epsilon-k(1-\eta-\epsilon ))n.\end{equation}

Let $$f(p)=\frac{r-3-p-\epsilon (C+r)}{r-3+\epsilon-p(1-\eta-\epsilon)}.$$  It is easy to see that $f(p)$ is decreasing with respect to $1\leq p\leq r-4$.
Then  (\ref{claim3-2}) and (\ref{claim3-3}) induce
\begin{eqnarray*}
  \frac{\sum\limits_{\substack{uv\in E(L\backslash A,S)\\ u\in L\backslash A}}x_u}{e(L\backslash A,S)} &\geq& f(k)\geq f(r-4)
     =\frac{1-\epsilon (C+r)}{1+\epsilon+(r-4)(\eta+\epsilon)}\\
     &\geq&1-(C+2r-3)(\eta+\epsilon).
\end{eqnarray*}
Hence there exists a vertex $v_{k+1}\in L\backslash U$ such that $$ x_{k+1}\geq1-(C+2r-3)(\eta+\epsilon).$$
Therefore by Claim~2,
\begin{eqnarray*}
d_{v_{k+1}} &\geq&(1- C_1((C+2r-3)(\eta+\epsilon)+\epsilon))n\geq (1- C_1(C+2r-2)(\eta+\epsilon))n.
\end{eqnarray*}

Let $c_{k+1}=C_1(c+2r-2)(c_2+1)$ which is independent of $\varepsilon$ and $n$. Then $x_i\ge 1-c_{k+1}\varepsilon$ and $d_{v_i}\ge (1-c_{k+1}\varepsilon)n$ for $i=1, \ldots, k+1$.
Hence by the induction principle, Claim~3 holds.



Denote by $A=\{v_1,v_2,\ldots,v_{r-3}\}$,  $B=\cap_{i=1}^{r-3} N_G(v_i)$ and $R=V(G)\backslash (A\cup B)$.
Then by $d_{v_i}\ge (1-C_3\varepsilon)n$ for $i=1, \ldots, r-3$, we have
 $$|B|\geq \sum_{i=1}^{r-3}(1-C_3\varepsilon)n-(r-4)n=(1-(r-3)C_3\epsilon)n$$
and
$$|R|=n-|A|-|B|\leq (r-3)C_3\epsilon n.$$


{\bf Claim~4.}
$G[A]$ induced by $A$ in $G$ is  a clique.

Clearly, $[A,B]$ is complete bipartite graph with $|A|=r-3$ and $|B|=(1-\delta)n$, where $\delta \le (r-3)C_3\epsilon$. Moreover,
$(1-3\delta)n>r$.
Since adding edges to a connected graph strictly increases its spectral radius, by Theorem~\ref{edge-minor-1}  and the maximality of $G$, $A$ must induce a  clique in $G$. This proves Claim~4.

{\bf Claim~5.}
For $v\in V(G)\backslash A$, $x_v\leq\frac{1}{C+3}$.

For any vertex $v\in R$, $|N_{G}(v)\cap B|\leq2$. Otherwise $A$, three edges $vu_i$ for $i=1,2,3$ form $K_r^-$ minor in ${G}$, where $u_i\in B$,  a contradiction. Furthermore, by definition of $R$, $|N_{G}(v)\cap A|\leq r-4$, which implies that
\begin{equation}\label{claim5-1}|N_{G}(v)\cap(A\cup B)|=|N_{G}(v)\cap B|+|N_{G}(v)\cap A|\leq r-4+2=r-2.\end{equation}
Hence
\begin{eqnarray*}
  \lambda \sum_{v\in R} x_v &=&  \sum_{v\in R} \sum_{uv\in E(G)} x_u
   \leq \sum_{v\in R}d_v\leq 2e(R)+e(R,A\cup B)\leq2e(R)+(r-2)|R|.
\end{eqnarray*}
Note that $G[R]$ is also $K_r^-$-minor free. Then
 \begin{equation}\label{claim5-2}
 \sum_{v\in R} x_v\leq\frac{2e(R)+(r-2)|R|}{\lambda}\leq \frac{C|R|+(r-2)|R|}{\lambda}=\frac{(C+r-2)|R|}{\lambda}.\end{equation}

For any vertex $v\in  B$, $|N_{G}(v)\cap B|\leq1$, otherwise $A$, $vw_1$ and $vw_2$ form  $K_r^-$ minor in $G$, where $w_1,w_2\in B$. This is a contradiction.  So  for any $v\in  B$,
\begin{equation}\label{claim5-3}
|N_{G}(v)\cap(A\cup B)|=|N_{G}(v)\cap B|+|N_{G}(v)\cap A|\leq 1+r-3=r-2.\end{equation}
Hence  by (\ref{claim5-1}), (\ref{claim5-3}) and (\ref{claim5-2}),  we have $|N_{G}(v)\cap(A\cup B)|\leq r-2$  for any $v\in V({G})\backslash A=R\cup B$  and
\begin{eqnarray*}
  \lambda x_v = \sum\limits_{uv\in E(G)} x_u\leq\sum_{\substack{uv\in E(G)\\u\in A\cup B}}x_u+\sum_{\substack{uv\in E(G)\\u\in R}}x_u \leq r-2+\sum_{u\in R}x_u
   \leq r-2+  \frac{(C+r-2)|R|}{\lambda},
\end{eqnarray*}
which implies that
\begin{eqnarray*}
  x_v &\leq& \frac{r-2}{\lambda}+ \frac{(C+r-2)|R|}{\lambda^2}
   \leq\frac{r-2}{\sqrt{(r-3)(n-r+3)}} +\frac{(C+r-2)C_3\epsilon n}{(r-3)(n-r+3)}\\
   &\leq&\frac{1}{2(C+3)}+\frac{1}{2(C+3)}
   =\frac{1}{C+3},
\end{eqnarray*}
where the last second inequality holds as  long as $n\geq \frac{4(r-2)^2(C+3)^2+r-3}{r-3}$ and $\epsilon\leq\frac{r-3}{2(C+3)(C+r-2)C_3}-\frac{(r-3)^2}{n}$.
So Claim~5 holds.


{\bf Claim~6.} $d_v=n-1$ for any $v\in A$.

Suppose that $d_v<n-1$. Since each vertex in $A$ is adjacent all vertices in $B$ by the definition of $A$ and $B$, there exists a vertex $z\in R$ such that $v\in A$ is not adjacent to $z\in R$.
 Let
  $$G^*=G-\{zw: w\in N_{B}(z)\}+\{zw :w\in A\setminus  N_A(z)  \}.$$
Not that $|N_G(z)\bigcap B|\le 2$. Hence  by Rayleigh Quotient and  Claim~5,
\begin{eqnarray*}
  \lambda(G^*)-\lambda(G) &\geq& \frac{{\bf x}^{\mathrm{T}}A(G^*){\bf x}}{{\bf x}^{\mathrm{T}}{\bf x}}-\frac{{\bf x}^{\mathrm{T}} A(G) {\bf x}}{{\bf x}^{\mathrm{T}}{\bf x}}
  \geq\frac{2x_z}{{\bf x}^{\mathrm{T}}{\bf x}}\bigg(x_v-\sum\limits_{\substack{zw\in E(G)\\w\in B}}x_w\bigg) \\
   &\geq&\frac{2x_z}{{\bf x}^{\mathrm{T}}{\bf x}}\bigg(x_v-\frac{2}{C+3}\bigg)\geq\frac{2x_z}{{\bf x}^{\mathrm{T}}{\bf x}}\bigg(1-C_3\epsilon-\frac{2}{C+3}\bigg)>0,
\end{eqnarray*}
where the last inequality holds as long as $\epsilon<\frac{1}{(C+3)C_3}$ .
Then $G^*$ is a $K_r^-$-minor free graph with larger spectral radius, which is a contradiction. Hence Claim~6 holds.

Now we can choose $\epsilon$ such that
$$ (\max\{\sqrt{C'C},C\})^{-\frac{1}{4}}(r-3.5)^{-\frac{1}{4}}n^{-\frac{1}{4}}\}\leq \epsilon\leq\frac{r-3}{2(C+3)(C+r-2)C_3}-\frac{(r-3)^2}{n}$$
  for sufficiently large $n$.

{\bf Claim~7.} The graph $G[B]$ induced by $B$ consists of independent edges  and isolated vertices.

By Claim~4, $G[A]$ is a induces.
If there is a path $P_3$ in $G[B]$, then  $A$ and this $P_3$ would form a $K_r^-$ minor in $G$, which is a contradiction.  So $G[B]$ does not contain $P_3$ and Claim~7 holds.

It follows from  Claims~4, 6 and 7 that  $G\subseteq F_{n,r-3}$. Note that adding an edge to  a connected graph will strictly increase spectral radius. By the maximality of $G$, $G$ must be $F_{n,r-3}$.

{\bf Case~2.} $l\geq2$. Note that $F$ is a subgraph of $K_r^-$. If $G$ is $F$-minor free, then $G$ is also $K_r^-$-minor free.
By Case~1, $G=F_{n,r-3}$. Since $F_{n,r-3}$ is also $F$-minor free, So the assertion holds.
\QEDB

\section{Proof of Theorem~\ref{thm1}~(ii) }

In order to prove Theorem~\ref{thm1}~(ii), we first prove the following two results.
\begin{Lemma}\label{L-minor-1}
Let  $G$ be a $K_r^{=}$-minor free graph of order $n$ with $ r\geq4$. If $G$ contains a complete bipartite subgraph $G[A,B]$  with $|A|=r-3$ and $|B|=(1-\delta)n\geq r$, then $B$ is an independent set and there are at least  $(1-2\delta)n$ vertices in $B$ which has no neighbors in $V(G)\backslash (A\cup B)$.
\end{Lemma}

\begin{Proof}
Suppose that there is an edge $uv$ in $G[B]$.  We can choose $w\in B\setminus\{u,v\}$ and $r-3$ edges $u_iv_i$ in $G$  with $u_i\in A$ and $v_i\in B\setminus \{u,v,w\}$ for $1\le i\le r-3$. Then $u_1v_1, \ldots, u_{r-3}v_{r-3}$, $uv$ and $w$ form a $K_r^{=}$-minor, which is a contradiction. Hence  $B$ is an independent set. Furthermore.

{\bf Claim.} If $C$ is a  component of $G-(A\cup B)$, then there is at most  a vertex in $B$ with a neighbor in $C$.

Suppose that there are two vertices $u$ and $v$ in $B$  with a neighbor in $C$. Then there is a path $P$ from $u$ to $v$ with all interior vertices in $C$. Let $w\in B\backslash \{u,v\}$.  Then $u_1v_1, \ldots, u_{r-3}v_{r-3}$, $P$ and $w$ form a
 $K_r^{=}$-minor, which is a contradiction. So Claim holds.

%
%
%
%

Denote by $R=V(G)\backslash (A\cup B)$ and
 $D=\{v\in B: N_G(v)\cap R=\emptyset \}.$
 By the definition of $R$, $$|R|=n-|A|-|B|\leq n-(r-3)-(1-\delta)n<\delta n,$$ which implies that
$R$ has at most $ \delta n$ components. By Claim, $B\backslash D$ has at most $\delta n$ vertices.
Hence $$|D|=|B|- |B\backslash D|\geq (1- \delta) n-\delta n\geq(1- 2\delta) n.$$
This completes the proof.
\end{Proof}

\begin{Theorem}\label{edge-minor-2}
Let  $G$ be a $K_r^{=}$-minor free graph of order $n$ with $ r\geq4$.  If there exist two disjoint subsets $A,B\subseteq V(G)$ such that $G[A]$  is a complete graph and $G[A,B]$ is a complete bipartite graph, where $|A|=r-3$, $|B|\geq1$,  then
 $G^*=G-\{uv:v\in N_B(u)\}+\{uv:v\in A\setminus N_A(u)\}$  is $K_r^{=}$-minor free for any  $u\in V(G)\backslash  (A\cup B)$, where $N_A(u)=N_G(u)\bigcap A$ and $N_B(u)=N_G(u)\bigcap B$

 \end{Theorem}

\begin{Proof}
Suppose that  $G^*$ has a  $K_r^{=}$-minor. Then there are $r$ disjoint  vertex subsets $X_1,\ldots, X_r\subseteq V(G^*)$ with the following properties:

(a) For each $X_i$, $G^*[X_i]$ is connected.

  (b)  
   There is at least one edge of $G^*$ between $X_i$ and $X_j$ for all distinct  $X_i$ and $X_j$ except two pairs $\{X_p ,X_s\}$ and $\{X_p,X_t\}$.


  Denote $D=\{v\in B: N_G(v)\cap (V(G)\backslash (A\cup B))=\emptyset \}$ and  $$f(X_1,\ldots,X_r)=|\{i:X_i\cap D\neq \emptyset~\text{for}~1\leq i\leq r\}|.$$
  We can choose  $X_1, \ldots, X_r$ such that  $f(X_1,\ldots,X_r)$ is as large as possible.
By Lemma~\ref{L-minor-1}, $D$ is an independent set with at least $(1-2\delta)n$ vertices.
By the similar to the proof of Theorem~\ref{edge-minor-1}, we can also choose  $X_1, \ldots, X_r$ such that  $f(X_1,\ldots,X_r)$ is as large as possible and  $|X_i\cap D|\leq1$ for all $1\leq i\leq r$.

Furthermore, by similar to the proof of Theorem~\ref{edge-minor-1}, $X_1,\ldots, X_r$ form a $K_r^{=}$-minor of $G$, which is a contradiction.  This completes the proof.
\end{Proof}

\begin{Theorem}\label{edge-minor-4}
Let  $G$ be a $K_r^{=}$-minor free graph of order $n$ with $ r\geq4$. If $G$ contains  a complete split subgraph  $S_{|A\cup B|,|A|}$ with partition $A$ and $B$, where $|A|=r-3$,  $|B|=(1-\delta)n$, and $(1-2\delta)n>r$.
If $G^*=G-\{uv:v\in N_G(u)\}+\{uv:v\in A\}$ for any $u\in V(G)\backslash (A\cup B)$,
then $G^*$ is also $K_r^{=}$-minor free.
 \end{Theorem}

 \begin{Proof} The proof of Theorem~\ref{edge-minor-3} is similar to Theorem~\ref{edge-minor-2} and omitted.\end{Proof}

Now we are ready to prove Theorem~\ref{thm1}~(ii).

{\bf Proof of Theorem~\ref{thm1}~(ii).}

 Let $G$ be an $F$-minor free graph of  sufficiently large order  $n$ with the maximum spectral radius.

{\bf Case~1.}
 $G$ is connected. Now we consider the following two subcases.

{\bf Case~1.1.} $l=1$, i.e, $F=K_r^{=}$.
let $\mathbf x=(x_u)_{u\in V(G}$ be a positive eigenvector with the maximum entry $x_w=1$ which corresponds  to   the  spectral radius  $\lambda (G)$.   For $0< \epsilon< 1$, let $$L=\{v\in V(G): x_v> \epsilon\}, S=\{v\in V(G): x_v\leq \epsilon\} . $$
Clearly,
$V(G)=L\cup S$.

By similar proofs of Claims~1-3 in the proof of  Theorem~\ref{thm1}~(i),
we can also choose $\epsilon$ small enough that there are   $r - 3$ vertices
 $A=\{v_1, \ldots, v_{r-3}\} \subseteq L$ which satisfy
$x_{v_i}\ge 1-C_3\varepsilon$ and $d_{v_i}\ge (1-C_3\varepsilon)n$ for $1\le i\le r-3$, where $C_3$ is  a constant independent of $\varepsilon$ and $n$.
Let $B=\bigcap_{i=1}^{r-3}N_G(v_i)$. By  similar proofs of Claims~4-6 in the proof of  Theorem~\ref{thm1}~(i),  we have  that $G[A]$ is a clique and  each vertex in $A$ has degree $n-1$, $B$ is an independent set.
Therefore  $G=S_{n,r-3}$.







{\bf Case~1.2.} $l\geq2$. Note that $F$ is a subgraph of $K_r^{=}$. If $G$ is $F$-minor free, then $G$ is also $K_r^{=}$-minor free.
By Case~1.1, $G=S_{n,r-3}$. Since $S_{n,r-3}$ is also $F$-minor free, the assertion holds.

{\bf Case~2.}
 $G$ is not connected. For $r\geq5$, adding an edge to  components of $G$ to make it connected  will result in a new $F$-minor free graph with larger spectral radius, which is a contradiction. For $r=4$, $F=K_1\nabla(K_1\cup K_2)$ or $F=P_4$ (in this situation, $F$ contains cut edges).
Let $G_1$ be a component of $G$  such that $\lambda(G_1)=\lambda(G)$. Set $n_1=|V(G_1)|$. Since $S_{n,r-3}$  is $F$ free, we have $\lambda(G)\geq\lambda(S_{n,r-3})>\sqrt{(r-3)(n-r+3)}$.
Then
\begin{eqnarray*}
  n_1-1 &\geq&  \lambda(G_1)=\lambda(G)\geq \sqrt{(r-3)(n-r+3)}
\end{eqnarray*}
which implies that
$n_1$ is also sufficiently large. By Case~1, $G_1=S_{n_1,r-3}$. Then
\begin{eqnarray*}
 \lambda(G) =\lambda(G_1) =\lambda(S_{n_1,r-3}) <\lambda(S_{n,r-3}),
\end{eqnarray*}
which is a contradiction. Therefore we finish the proof of  Theorem~\ref{thm1}~(ii).
\QEDB

\section{Conclusion and discussion}
Indeed,  Theorem~\ref{thm1} can be extended as follows.
\begin{Theorem}\label{thm}
Let  $r\geq 4$ and $M \subseteq E(K_r)$ such that $G$ is a ($K_r- M$)-minor free graph of  order  $n$. \\
(i) If $M$ is a matching, then  $\lambda(G)\leq \lambda(F_{n,r-3})$ with equality if and only if $G=F_{n,r-3}$ for sufficiently large $n$.\\
(ii) $M$ is not a matching and $G[M]$  is connected  containing no $K_3$, then  $\lambda(G)\leq  \lambda(S_{n,r-3})$ with equality if and only if $G=S_{n,r-3} $ for sufficiently large $n$.
\end{Theorem}

\begin{Proof}
Let $G$ be a ($K_r- M$)-minor free graph of  sufficiently large order  $n$ with the maximum spectral radius.

(i) Since  $K_r- M$ is a subgraph of $K_r^-$, $G$ is also  $K_r^-$-minor free. By Theorem~\ref{thm1}~(i), $G=F_{n,r-3}$. Note that $F_{n,r-3}$ is also ($K_r- M$)-minor free, the assertion holds.

(ii) Since  $K_r-M$ is a subgraph of $K_r^{=}$, $G$ is also  $K_r^{=}$-minor free. By Theorem~\ref{thm1}~(ii), $G=S_{n,r-3}$. Note that $S_{n,r-3}$ is also ($K_r-M$)-minor free, the assertion holds.
\end{Proof}

By Theorem~\ref{thm}~(ii), we have the following two corollaries.
\begin{Corollary}
 If $G$ is a $(K_{r}-E(C_r))$-minor free graph of  sufficiently large order  $n$ with $r\geq 4$.,
then  $\lambda(G)\leq  \lambda(S_{n,r-3})$ with equality if and only if $G=S_{n,r-3}$.
\end{Corollary}

\begin{Corollary}\label{cor}\cite{Tait2019}
 If $G$ is a $K_{r}$-minor free graph of  sufficiently large order  $n$ with $r\geq 3$,
then  $\lambda(G)\leq  \lambda(S_{n,r-2})$ with equality if and only if $G=S_{n,r-2}$.
\end{Corollary}

\begin{Proof}
Let $G$ be a $K_{r}$-minor free graph of  sufficiently large order  $n$ with the maximum spectral radius.
Since $G$ is $K_{r}$-minor, $G$ is also  ($K_{r}\cup K_1)$-minor free. Note that ($K_{r}\cup K_1$) is obtained from $K_{r+1}$ by deleting $r$ edges with a common vertex.
By  Theorem~\ref{thm1}~(ii), $G=S_{n,r-3}$. Note that $S_{n,r-3}$ is also $K_{r}$-minor free, the result follows directly.
 \end{Proof}

  \vspace{2mm}
Theorem~\ref{thm} is very interesting. Next we discuss Theorem~\ref{thm} for $r=4$ and $r=5$.

{\bf Remark~1.} (i) Note that $|M|\leq2$. If $|M|=1$, then $K_4-M=F_1$ (see Fig.1).   If $|M|= 2$, then $K_4-M=F_2$. By Theorem~\ref{thm}~(i),
$F_{n,1}$ is the graph with the maximum  spectral radius among  $F_i$-minor free graphs of sufficiently large order $n$ for $i=1,2$.

(ii) Note that $2\leq|M|\leq4$. If $|M|=2$, then $K_4-M=F_3$. If  $|M|=3$, then  $K_4-M\in \{F_4,F_5\}$. If  $|M|=4$, then  $K_4-M=F_6$.
By Theorem~\ref{thm},
$S_{n,1}$ is the graph with the maximum  spectral radius among $F_i$-minor free graphs of sufficiently large order $n$ for $i=3,4,5,6$.

  \vspace{2mm}
{\bf Remark~2.} (i)  Note that $|M|\leq2$. If $|M|=1$, then $K_5-M=F_7$.  If $|M|= 2$, then $K_4-M=F_8$. By Theorem~\ref{thm}~(i),
$F_{n,2}$ is the graph with the maximum  spectral radius among $F_i$-minor free graphs of sufficiently large order $n$ for $i=7,8$.

(ii)  Note that $2\leq|M|\leq6$. If $|M|=2$, then $K_5-M=F_9$. If  $|M|=3$, then  $K_5-M\in\{F_{10},F_{11},F_{12}\}$. If  $|M|=4$, then  $K_5-M\in\{F_{13},F_{14},F_{15},F_{16}\}$. If  $|M|=5$, then  $K_5-M\in\{F_{17},F_{18}\}$.  If  $|M|=6$, then  $K_5-M=F_{19}$.
By Theorem~\ref{thm}~(ii),
$S_{n,3}$ is the graph with the maximum  spectral radius among $F_i$-minor free graphs of sufficiently large order $n$ for $i=9,\dots,19$.

\vskip 1cm

\begin{center}
\begin{tikzpicture}[scale=1]

\draw((0,0)--(1,0)--(0.5,0.65)--(0,0)--(0.5,-0.65)--(1,0);

\draw[fill](0,0)circle [radius=0.05];
\draw[fill](1,0)circle [radius=0.05];
\draw[fill](0.5,0.65)circle [radius=0.05];
\draw[fill](0.5,-0.65)circle [radius=0.05];

\draw(0.5,-1.2)node{\scriptsize  $F_1$};

\draw((2,0)--(2.5,0.65)--(3,0)--(2.5,-0.65)--(2,0);

\draw[fill](2,0)circle [radius=0.05];
\draw[fill](3,0)circle [radius=0.05];
\draw[fill](2.5,0.65)circle [radius=0.05];
\draw[fill](2.5,-0.65)circle [radius=0.05];

\draw(2.5,-1.2)node{\scriptsize  $F_2$};

\draw((4,0)--(5,0)--(4.5,0.65)--(5,0)--(4.5,-0.65)--(4,0);

\draw[fill](4,0)circle [radius=0.05];
\draw[fill](5,0)circle [radius=0.05];
\draw[fill](4.5,0.65)circle [radius=0.05];
\draw[fill](4.5,-0.65)circle [radius=0.05];

\draw(4.5,-1.25)node{\scriptsize  $F_3$};

\draw(6,0)--(6.5,-0.65)--(7,0)--(6,0);

\draw[fill](6,0)circle [radius=0.05];
\draw[fill](7,0)circle [radius=0.05];
\draw[fill](6.5,0.65)circle [radius=0.05];
\draw[fill](6.5,-0.65)circle [radius=0.05];

\draw(6.5,-1.2)node{\scriptsize  $F_4$};

\draw(8.5,0.65)--(8,0)--(8.5,-0.65)--(9,0);

\draw[fill](8,0)circle [radius=0.05];
\draw[fill](9,0)circle [radius=0.05];
\draw[fill](8.5,0.65)circle [radius=0.05];
\draw[fill](8.5,-0.65)circle [radius=0.05];

\draw(8.5,-1.2)node{\scriptsize  $F_5$};

\draw(10,0.45)--(11,0.45);

\draw(10,-0.45)--(11,-0.45);

\draw[fill](10,0.45)circle [radius=0.05];
\draw[fill](11,0.45)circle [radius=0.05];
\draw[fill](10,-0.45)circle [radius=0.05];
\draw[fill](11,-0.45)circle [radius=0.05];

\draw(10.5,-1.2)node{\scriptsize  $F_6$};

\draw(12,0)--(13,0)--(12.5,0.5)--(12.75,-0.55)--(12.25,-0.55)--(12.5,0.5);
\draw(12,0)--(12.75,-0.55)--(13,0)--(12.25,-0.55)--(12,0);

\draw[fill](12.5,0.5)circle [radius=0.05];
\draw[fill](12,0)circle [radius=0.05];
\draw[fill](13,0)circle [radius=0.05];
\draw[fill](12.75,-0.55)circle [radius=0.05];
\draw[fill](12.25,-0.55)circle [radius=0.05];

\draw(12.5,-1.2)node{\scriptsize  $F_7$};

\draw(0,-2.25)--(1,-2.25)--(1,-3.25)--(0,-3.25)--(0,-2.25)--(0.5,-2.75)--(1,-2.25)--(0.5,-2.75)--(1,-3.25)--(0.5,-2.75)--(0,-3.25);
\draw[fill](0,-2.25)circle [radius=0.05];
\draw[fill](1,-2.25)circle [radius=0.05];
\draw[fill](1,-3.25)circle [radius=0.05];
\draw[fill](0,-3.25)circle [radius=0.05];
\draw[fill](0.5,-2.75)circle [radius=0.05];

\draw(0.5,-3.65)node{\scriptsize  $F_8$};

\draw(3,-2.25)--(3,-3.25)--(2,-3.25)--(2,-2.25)--(3,-2.25)--(2,-3.25)--(3,-3.25)--(2,-2.25);
\draw(3,-2.25)--(3.5,-2.75)--(3,-3.25);

\draw[fill](2,-2.25)circle [radius=0.05];
\draw[fill](3,-2.25)circle [radius=0.05];
\draw[fill](3,-3.25)circle [radius=0.05];
\draw[fill](2,-3.25)circle [radius=0.05];
\draw[fill](3.5,-2.75)circle [radius=0.05];

\draw(2.5,-3.65)node{\scriptsize  $F_9$};

\draw(4.5,-2.25)--(5.5,-3.25)--(4.5,-3.25)--(4.5,-2.25)--(5.5,-2.25)--(4.5,-3.25)--(5.5,-3.25)--(4.5,-2.25);
\draw(6,-2.75)--(5.5,-2.25)--(5.5,-3.25);

\draw[fill](4.5,-2.25)circle [radius=0.05];
\draw[fill](5.5,-2.25)circle [radius=0.05];
\draw[fill](5.5,-3.25)circle [radius=0.05];
\draw[fill](4.5,-3.25)circle [radius=0.05];
\draw[fill](6,-2.75)circle [radius=0.05];
\draw(5,-3.65)node{\scriptsize  $F_{10}$};

\draw(8,-3.25)--(7,-3.25)--(7,-2.25)--(8,-2.25)--(7,-3.25)--(8,-3.25)--(7,-2.25);
\draw(8,-2.25)--(8.5,-2.75)--(8,-3.25);

\draw[fill](7,-2.25)circle [radius=0.05];
\draw[fill](8,-2.25)circle [radius=0.05];
\draw[fill](8,-3.25)circle [radius=0.05];
\draw[fill](7,-3.25)circle [radius=0.05];
\draw[fill](8.5,-2.75)circle [radius=0.05];

\draw(7.75,-3.65)node{\scriptsize  $F_{11}$};

\draw(10.5,-3.25)--(9.5,-3.25)--(9.5,-2.25)--(10.5,-2.25)--(10.5,-3.25)--(9.5,-2.25);
\draw(10.5,-3.25)--(10.5,-2.25)--(11,-2.75)--(10.5,-3.25);

\draw[fill](9.5,-2.25)circle [radius=0.05];
\draw[fill](10.5,-2.25)circle [radius=0.05];
\draw[fill](10.5,-3.25)circle [radius=0.05];
\draw[fill](9.5,-3.25)circle [radius=0.05];
\draw[fill](11,-2.75)circle [radius=0.05];

\draw(10,-3.65)node{\scriptsize  $F_{12}$};

\draw(13,-3.25)--(12,-3.25)--(12,-2.25)--(13,-2.25)--(13,-3.25)--(12,-2.25);
\draw(13,-2.25)--(12,-3.25);

\draw[fill](12,-2.25)circle [radius=0.05];
\draw[fill](13,-2.25)circle [radius=0.05];
\draw[fill](13,-3.25)circle [radius=0.05];
\draw[fill](12,-3.25)circle [radius=0.05];
\draw[fill](13.5,-2.75)circle [radius=0.05];

\draw(12.5,-3.65)node{\scriptsize  $F_{13}$};

\draw(1,-5.75)--(0,-4.75)--(1,-4.75)--(1,-5.75)--(0,-5.75)--(0,-4.75);
\draw(1,-4.75)--(1.5,-5.25);

\draw[fill](0,-4.75)circle [radius=0.05];
\draw[fill](1,-4.75)circle [radius=0.05];
\draw[fill](1,-5.75)circle [radius=0.05];
\draw[fill](0,-5.75)circle [radius=0.05];
\draw[fill](1.5,-5.25)circle [radius=0.05];
\draw(0.5,-6.15)node{\scriptsize  $F_{14}$};

\draw(2.5,-4.75)--(3.5,-4.75)--(3.5,-5.75)--(2.5,-5.75)--(2.5,-4.75);
\draw(3.5,-4.75)--(4,-5.25)--(3.5,-5.75);

\draw[fill](2.5,-4.75)circle [radius=0.05];
\draw[fill](3.5,-4.75)circle [radius=0.05];
\draw[fill](3.5,-5.75)circle [radius=0.05];
\draw[fill](2.5,-5.75)circle [radius=0.05];
\draw[fill](4,-5.25)circle [radius=0.05];

\draw(3,-6.15)node{\scriptsize  $F_{15}$};

\draw(5,-4.75)--(5.75,-5.25)--(6.5,-4.75)--(6.5,-5.75)--(5.75,-5.25)--(5,-5.75)--(5,-4.75);
\draw[fill](5,-4.75)circle [radius=0.05];
\draw[fill](5.75,-5.25)circle [radius=0.05];
\draw[fill](6.5,-5.75)circle [radius=0.05];
\draw[fill](5,-5.75)circle [radius=0.05];
\draw[fill](6.5,-4.75)circle [radius=0.05];
\draw(5.75,-6.15)node{\scriptsize  $F_{16}$};

\draw(7.5,-4.75)--(7.5,-5.75)--(8,-5.25)--(7.5,-4.75);
\draw(8,-5.25)--(8.5,-5.75)--(8.5,-4.75);

\draw[fill](7.5,-4.75)circle [radius=0.05];
\draw[fill](8,-5.25)circle [radius=0.05];
\draw[fill](7.5,-5.75)circle [radius=0.05];
\draw[fill](8.5,-5.75)circle [radius=0.05];
\draw[fill](8.5,-4.75)circle [radius=0.05];
\draw(8,-6.15)node{\scriptsize  $F_{17}$};

\draw(9.5,-5)--(10.25,-4.5)--(11,-5)--(10.75,-5.75)--(9.75,-5.75)--(9.5,-5);
\draw[fill](9.5,-5)circle [radius=0.05];
\draw[fill](10.25,-4.5)circle [radius=0.05];
\draw[fill](11,-5)circle [radius=0.05];
\draw[fill](10.75,-5.75)circle [radius=0.05];
\draw[fill](9.75,-5.75)circle [radius=0.05];
\draw(10.25,-6.15)node{\scriptsize  $F_{18}$};

\draw(12,-4.75)--(12,-5.75);
\draw(12.5,-4.75)--(13,-5.25)--(12.5,-5.75)--(12.5,-4.75);
\draw[fill](12.5,-4.75)circle [radius=0.05];
\draw[fill](12,-5.75)circle [radius=0.05];
\draw[fill](13,-5.25)circle [radius=0.05];
\draw[fill](12.5,-5.75)circle [radius=0.05];
\draw[fill](12,-4.75)circle [radius=0.05];
\draw(12.5,-6.15)node{\scriptsize $F_{19}$};

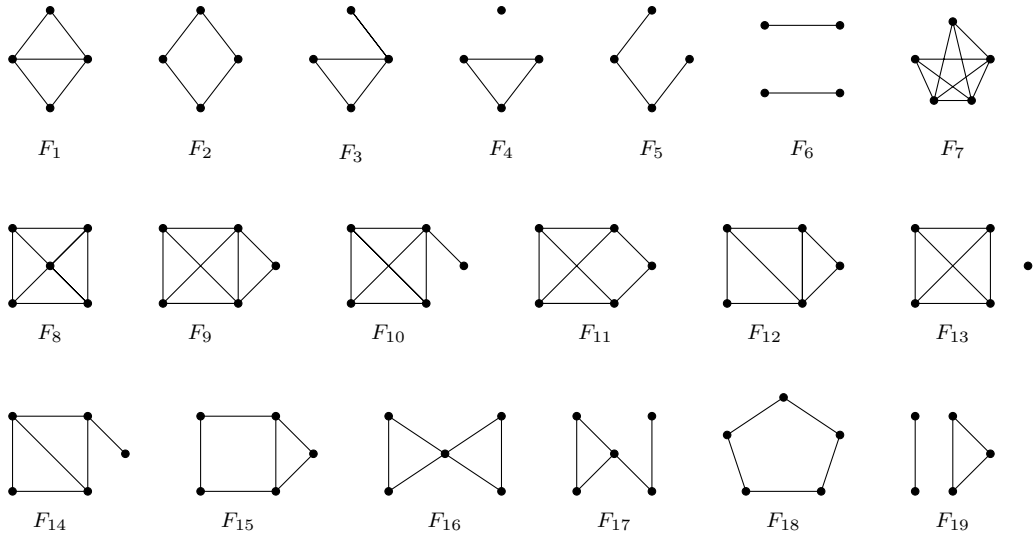
\draw(6,-7.2)node{ Fig.1. Graphs $F_1$--$F_{19}$};
\end{tikzpicture}
\end{center}

\end{document}